\documentstyle[12pt]{article}
\def\R{\ifmmode{\rm I\mkern-3.1mu
R\mkern1mu}\else{\rm I\kern-.18em  R\hskip1pt\ 
}\fi\relax}  
 
\def\b{\beta} 

\def\G{\Gamma}

\def\t{\tau}
\def\d{\delta}  
\def\th{\theta} 
\def\l{\lambda}

\def\n{\nu}
\def\di{\diamond}

\def\s{\sigma}
\def\sou{\overline}
\def\so{\underline} 
\def\O{\Omega}
\def\f{\rightarrow}
\def\q{\forall}
\def\e{\exists}  
\def\v{\vdash}
 
\def\p{\succ}
 
\def\mats{\ifmmode{ {\hbox{\bigreek s}} }\else{ 
{\bigreek s} }\fi\relax}
\def\matsin{\ifmmode{ {\hbox{\smgreek s}} }\else{ 
{\smgreek s} }\fi\relax}
\def\matt{\ifmmode{ {\hbox{\bigreek t}} }\else{ 
{\bigreek t} }\fi\relax}
\def\mattin{\ifmmode{ {\hbox{\smgreek t}} }\else{ 
{\smgreek t} }\fi\relax}

\newtheorem{theo}{Theorem}[section]
\newtheorem{lemma}{Lemma}[section]

\begin{document} 

\vspace*{0cm}
\begin{center} 
\LARGE\bf 
A General Type for Storage Operators\\[1cm]
\end{center}

\begin{center} 
\bf 
Karim NOUR\\
\rm
LAMA - Equipe de Logique, Universit\'e de Chamb\'ery\\
73376 Le Bourget du Lac\\
e-mail nour@univ-savoie.fr\\[1cm]
\end{center}

\begin{abstract}
In 1990, J.L. Krivine introduced the notion of storage operator to simulate, in $\l$-calculus, the
"call by value" in a context of a "call by name". J.L. Krivine has shown that, using G\H{o}del
translation from classical into intuitionistic logic, we can find a simple type for storage
operators in $AF2$ type system.\\ 
In this present paper, we give a general type for storage operators in a slight extension of $AF2$. 
We give at the end (without proof) a generalization of this result to other types.
\end{abstract}

\section{Introduction}

In 1990, J.L. Krivine introduced the notion of storage operators (see [3]). They are closed $\l$-terms
which allow, for a given data type (the type of integers, for example), to simulate in $\l$-calculus
the "call by value" in a context of a "call by name" (the head reduction).\\

J.L. Krivine has shown that the formula $\q x \{ N$*$[x] \f \neg\neg N[x] \}$ is a specification for
storage operators for Church integers : where $N[x]$ is the type of integers in second order
logic, and the operation $*$ is the simple  G\H{o}del translation from classical into intuitionistic
logic which associates to every formula $F$ the formula $F$* obtained by replacing in $F$ each atomic
formula with its negation (see [3]).\\

Some authors have been interested in the research of a most general type for storage operators. For
example, V. Danos and L. Regnier have given as type for storage operators the formula $\q x \{N^e [x]
\f \neg\neg N[x] \}$ where the operation $e$ is an elaborate G\H{o}del translation which associates
to every formula $F$ the formula $F^e$ obtained by replacing in $F$ each atomic formula
$X(\sou{t})$ by $X_1(\sou{t}),...,X_r(\sou{t}) \f \perp$ (see [1]). J.L.
Krivine and the author have given a more general type for storage operators : the formula $\q x \{ N^g
[x] \f \neg\neg N[x] \}$ where the operation $g$ is the general G\H{o}del translation which
associates to every formula $F$ the formula $F^g$ obtained by replacing in $F$ each atomic formula
$X(\sou{t})$ by a formula $G_X [ \sou{t} / \sou{x} ]$ ending with $\perp$ (see [4]
and [5]).\\

With the types cited before, we cannot type the following simple storage operators :\\ 
$T=\l \n \l f(( \n )(T_i) \n f) \l xx$ and $T'=\l \n \l f (( \n )(T_i) \n f)\l d(T_j) \n f$ where
$T_i$ ($i=1$ or $2$) are the standard storage operators for integers (see [3]). This is due to the
fact that the normal form of $T$ (and $T'$) contains a variable $\nu$ applied to two arguments and
another $\nu$ applied to three arguments. Therefore, we cannot type $T$ and $T'$ because the variable
$\nu$ is assigned by $N$*$[x]$ (for example) and thus the number of the $\nu$-arguments is fixed once
for all.\\

To solve the problem, we will replace $N$*$[x]$ in the type of storage operators by another type
$N^{\perp}[x]$ which does not limit the number of $\nu$-arguments and only enables to generate
formulas ending with $\perp$ in order to find a general specification for storage operators.\\ 

The specifications of storage operators that we have obtained up to now do not explain that these
operators only accept integers (for example $\l n\l f \l z (x)z$ is a normal $\l$-term of type
$N$*$[0]$). We will see that the type $N^{\perp}[x]$ is also a specification for the integers.\\

In this paper, we give a general type for the storage operators for integers in a slight extension
of $AF2$ (the storage operators $T$ and $T'$ are typable of this type). We give at the end (without
proof) a generalization of this result to the $\q$-positive types (the universal second order
quantifier appears positively in these types). \\

\bf Acknowledgement.\/ \rm We wish to thank J.L. Krivine for helpful discussions. He 
found independently the principal result of this paper which he proved by a semantical method.
\\[3cm]

\section{Definitions and notations}

\subsection{The pure $\l$-calculus}
 
Let $t,u,u_1,...,u_n$ be $\l$-terms, the application of $t$ to $u$ is denoted by $(t)u$. In the same
way we write $(t)u_1...u_n$ instead of $(...((t)u_1)...)u_n$.\\  
The $\b$-equivalence relation is denoted by $u \simeq\sb{\b} v$.\\ 
The notation $\s(t)$ represents the result of the simultaneous substitution $\s$ to the free variables
of $t$ after a suitable renaming of the bounded variables of $t$.\\    
We denote by $(u)^n v$ the $\l$-term $(u)...(u)v$ where $u$ occurs $n$ times, and $\sou{u}$ the
sequence of $\l$-terms $u_1,...,u_n$. If $\sou{u} = u_1,...,u_n$, we denote by $(t)\sou{u}$ the
$\l$-term $(t)u_1...u_n$. \\  
Let us recall that a $\l$-term $t$ either has a head redex [i.e. $t=\l x_1 ...\l
x_n (\l x u) v v_1 ... v_m$, the head redex being $(\l x u) v$], or is in head normal form [i.e.
$t=\l x_1 ...\l x_n (x) v_1 ... v_m$].\\  
The notation $u \p v$ means that $v$ is obtained from $u$ by some head reductions.\\ 
If $u \p v$, we denote by $h(u,v)$ the length of the head reduction between $u$ and $v$. \\ 
A $\l$-term $t$ is said solvable iff the head reduction of $t$ terminates.

\begin{lemma} (see [3]) If $u \p v$, then :\\
1) for every substitution $\s$, $\s(u) \p \s(v)$ and $h(\s(u),\s(v))=h(u,v)$. \\
2) for every sequence of $\l$-terms $\sou{w}$, there is a $w$ such that $(u) \sou{w} \p w$, $(v)
\sou{w} \p w$, and $h((u) \sou{w},w)=h((v) \sou{w},w)+h(u,v)$.  
\end{lemma}

\subsection{The $AF2$ type system}

The types will be formulas of second order predicate logic over a given language.\\
The logical connectives are $\perp$ (for absurd), $\f$, and $\q$.\\
There are individual (or first order) variables denoted by $x,y,z,...,$ and predicate (or second
order) variables denoted by $X,Y,Z,....$\\ 
We do not suppose that the language has a special constant for equality. Instead, we define the
formula $u=v$ (where $u,v$ are terms) to be $\q Y(Y(u) \f Y(v))$ where $Y$ is a unary predicate
variable. Such a formula will be called an equation.\\
The formula $F_1 \f (F_2 \f(...\f (F_n \f G)...))$ is also denoted by $F_1,F_2,...,F_n \f G$.\\   
For every formula $A$, we denoted by $\neg A$ the formula $A \f \perp$.\\ 
If $\sou{v} = v_1,...,v_n$ is a sequence of variables, we denoted by $\q \sou{v} A$ the formula $\q
v_1...\q v_n A$.\\  
Let $t$ be a $\l$-term, $A$ a type, and $\G = x_1 : A_1 ,..., x_n : A_n$ a context. We
define by the mean of this following rules the notion " $t$ is of type $A$ in the context $\G$ ". This
notion is denoted by $\G\v t:A$.     
\begin{itemize} 
\item[] (1) $\G\v x_i:A_i$ $1 \leq i \leq n$. 
\item[] (2) If $\G,x:A \v t:B$, then $\G\v\l xt:A \f B$.  
\item[] (3) If $\G\v u:A \f B$, and $\G\v v:A$, then $\G\v (u)v:B$.
\item[] (4) If $\G\v t:A$, then $\G\v t: \q xA$. (*) 
\item[] (5) If $\G\v t:\q xA$, then $\G\v t:A[u/x]$. (**)   
\item[] (6) If $\G\v t:A$, then $\G\v t: \q XA$. (*) 
\item[] (7) If $\G\v t:\q XA$, then $\G\v t:A[G/X]$. (**) 
\item[] (8) If $\G\v t:A[u/x]$, then $\G\v t:A[v/x]$. (***) 
\end{itemize} 
The previous rules are subject to the following restrictions :\\
(*) The variable $x$ (resp. $X$) has no free occurence in $\G$.\\
(**) $u$ is a term and $G$ is a formula of the language.\\
(***) $u$ and $v$ are terms such that $u=v$ is a consequence of a given set of equations.\\

This type $\l$-calculus system is called $AF2$ (for arithm\'etique fonctionnelle du second ordre).
 
\begin{theo} (see [2]) The $AF2$ type system has the following properties :\\
1) Type is preserved during reduction.\\
2) Typable $\l$-terms are strongly normalizable.
\end{theo}

We define on the set of types the two binary relations $\lhd$ and $\approx$ as the least reflexive
and transitive binary relations such that :  
\begin{itemize} 
\item[] - $\q xA \lhd A[u/x]$, if $u$ is a term of language ;
\item[] - $\q XA \lhd A[F/X]$, if $F$ is a formula of language ;
\item[] - $A \approx B$ iff $A=C[u/x]$, $B=C[v/x]$, and $u=v$ is a consequence of a given set of
equations. 
\end{itemize}  

\begin{theo} (see [5] and [7]) \\
1) Let $A$ be an atomic formula. If $\G\v t:A$, then $t$ does not begin by
$\l$. \\ 2) If $\G ,x:A \v (x)u_1...u_n:B$, then :\\
$n=0$, $A \lhd C$, $C \approx C'$, $B=\q \sou{v}C'$, and $\sou{v}$ have no free occurence in $\G$ and
$A$,\\  
or\\  
$n \geq 1$, $ A \lhd C_1 \f B_1$, $B'_i \lhd C_{i+1} \f B_{i+1}$ $1 \leq i \leq n-1$, $B'_n \lhd
B_{n+1}$, $B= \q \sou{v} B'_{n+1}$ where $B_i \approx B'_i$ $1 \leq i \leq n+1$, $\G,x:A \v u_i:C_i$
$1 \leq i \leq n$, and $\sou{v}$ have no free occurence in $\G$ and $A$.  
\end{theo}

\section{The Church integers}

Each data type can be defined by a second order formula. For example, the type of integers is the
formula : 
\begin{center}
$N[x]= \q X \{ X(0), \q y(X(y) \f X(sy)) \f X(x) \}$ 
\end{center}
where $X$ is a unary predicate variable, $0$ is a constant symbol for zero, and $s$ is a unary
function symbol for successor.\\ The formula $N[x]$ means semantically that $x$ is an integer iff $x$
belongs to each set $X$ containing $0$ and closed under the successor function $s$.\\ 
The $\l$-term $\so{0} = \l x \l fx$ is of type $N[0]$ and represents zero.\\ 
The $\l$-term $\so{s} = \l n\l x\l
f(f)((n)x)f$ is of type $\q y(N[y] \f N[s(y)])$ and represents the successor function.\\   
A set of equations $E$ is said adequate with the type of integers iff :
\begin{itemize} 
\item[] - $s(a)=0$ is not an equational consequence of $E$ ;
\item[]- If $s(a)=s(b)$ is an equational consequence of $E$, then so is $a=b$.
\end{itemize}
In the rest of the paper, we assume that all the set of equations are adequate with the type of
integers.\\

For each integer $n$, we define the Church integer $\so{n}$ by $\so{n} = \l x\l f(f)^n x$.

\begin{theo} (see [2])
For each integer $n$, $\so{n}$ is the unique normal $\l$-term of type $N[s^n(0)]$.
\end{theo}

The propositional trace 
\begin{center}
$N=\q X \{ X,(X \f X) \f X \}$ 
\end{center}
of $N[x]$ also defines the integers.

\begin{theo} (see [2])
A normal $\l$-term is of type $N$ iff it is of the form $\so{n}$, for a certain integer $n$.
\end{theo} 

\bf Remark\/ \rm A very important property of data type is the following (we express it for the type
of integers) : in order to get a program for a function $f : N \f N$ it is sufficient to prove $\v
\q x ( N[x] \f N[f(x)] )$. For example a proof of $\v \q x ( N[x] \f N[p(x)] )$ from the equations
$p(0)=0$, $p(s(x))=x$ gives a $\l$-term for the predecessor in Church intergers (see [2]). $\Box$

\section{The storage operators}

Let $T$ be a closed $\l$-term. We say that $T$ is a storage operator for the integers iff for every
$n \geq 0$, there is $\t_n \simeq\sb{\b} \so{n}$, such that for every 
$\th_n \simeq\sb{\b} \so{n}$, there is a substitution $\s$, such that $(T)\th_n f \p (f)\s(\t_n)$.\\

\bf Remark\/ \rm 
Let $F$ be any $\l$-term (for a function), and $\th_n$ a $\l$-term $\b$-equivalent to $\so{n}$. During
the computation of $(F)\th_n$, $\th_n$ may be computed each time it comes in head position. Instead
of computing $(F)\th_n$, let us look at the head reduction of $(T)\th_n F$. Since it is
$\{(T)\th_n f\}[F/f]$, by Lemma 2.1, we shall first reduce $(T)\th_n f$ to its head normal form,
which is $(f)\s(\t_n)$, and then compute $(F)\s'(\t_n)$. The computation has been decomposed into two
parts, the first being independent of $F$. This first part is essentially a computation of $\th_n$,
the result being $\t_n$, which is a kind of normal form of $\th_n$. The substitutions made in $\t_n$
have no computational significance, since $\so{n}$ is closed. So, in the computation of $(T)\th_n F$,
$\th_n$ is computed first, and the result is given to $F$ as an argument, $T$ has stored the result,
before giving it, as many times as needed, to any function. $\Box$ \\

\bf Examples\/ \rm 
If we take :\\ $T_1 = \l n((n)\d)G$ where
$G = \l x\l y(x)\l z(y)(\so{s})z$ and $\d = \l f(f)\so{0}$ \\
$T_2 = \l n\l f(((n)f)F)\so{0}$
where $F = \l x\l y(x)(\so{s})y$, \\
then it is easy to check that : for every $\th_n \simeq\sb{\b} \so{n}$, $(T_i)\th_n f \p
(f)(\so{s})^n \so{0}$ ($i=1$ or $2$). \\ 
Therefore $T_1$ and $T_2$ are two storage operators for the integers. $\Box$ \\

It is a remarkable fact that we can give simple types to storage operators for integers. We first
define the simple G\H{o}del translation $F$* of a formula $F$ : it is obtained by replacing in the
formula $F$, each atomic formula $A$ by $\neg A$. For example :
\begin{center}
$N$*$[x]=\q X \{\neg X(0),\q y(\neg X(y) \f \neg X(sy)) \f \neg X(x) \}$
\end{center}
It is well know that, if $F$ is provable in classical logic, then $F$* is provable in intuitionistic
logic.\\ 

We can check that $\v T_1,T_2 : \q x \{N$*$[x] \f\neg\neg N[x] \}$. And, in general, we have the
following Theorem :

\begin{theo} (see [3] and [6])
If $\v T: \q x\{N$*$[x] \f\neg\neg N[x]\}$, then $T$ is a storage operator for the integers.
\end{theo}

\bf Remark\/ \rm Let $\th_0 = \l x \l f \l z (x) ( \l d z) \l xx$. \\
It is easy to check that $\v \th_0 : N$*$[0]$, and $( T_2 ) \th_0 f \p (f) ( \l d \so{0} ) \l xx$.\\
Therefore $T_2$ is not a storage operator for the set $\{ t$ / $ \v t : N$*$[s^n(0)]$ $n \geq 0
\}$. $\Box$ \\ 

The previous definition is not well adapted to study the storage operators. Indeed, it is, a priori, a
$\Pi^{0}_{4}$ statement ($\q n\e \t_n\q \th_n\e \s A(T,n,\t_n,\th_n,\s)$). We will show (Theorem 4.2)
that it is in fact equivalent to a $\Pi^{0}_{1}$ statement ($\t_n$ can be computed from $n$, and $\s$
from $\th_n$).\\

Let $\n$ and $f$ two fixed variables.\\ 
We denoted by $x_{n,a,b,\sou{c}}$ (where $n$ is an integer, $a,b$ two $\l$-terms, and $\sou{c}$ a
finite sequence of $\l$-terms) a variable which does not appear in $a,b,\sou{c}$.

\begin{theo} (see [5] and [8])
A closed $\l$-term $T$ is a storage operators for the integers iff for every $n \geq 0$, there is a
finite sequence of head reduction $\{ U_i \p V_i \}_{1\leq i\leq r}$ such that :\\ 
1) $U_1 = (T)\n f$ and $V_r = (f)\t_n$ where $\t_n \simeq\sb{\b}\so{n}$ ;\\
2) $V_i = (\n) a b \sou{c}$ or $V_i = (x_{l,a,b,\sou{c}}) \sou{d}$  $0 \leq l \leq n-1$;\\
3) If $V_i = (\n)a b \sou{c}$, then $U_{i+1} = (a)\sou{c}$ if $n=0$ and
$U_{i+1} = ((b)x_{n-1,a,b,\sou{c}})\sou{c}$ if $n \neq 0$ ;\\ 
4) If $V_i = (x_{l,a,b,\sou{c}})\sou{d}$ $0 \leq l \leq n-1$, then $U_{i+1} = (a)\sou{d}$ if $l=0$ and
$U_{i+1} = ((b)x_{l-1,a,b,\sou{d}})\sou{d}$ if $l \neq 0$.
\end{theo}

\section{General type for storage operators}

\subsection{The $AF2_{\perp}$ type system}

In this section, we present a slight extension of the $AF2$ type system denoted by $AF2_{\perp}$.\\

We assume that for every integer $n$, there is a countable set of special $n$-ary second order
variables denoted by $X_{\perp},Y_{\perp},Z_{\perp}$...., and called $\perp$-variables.\\ 

A type $A$ is called an $\perp$-type iff $A$ is obtained by the following rules : 
\begin{itemize} 
\item[] - $\perp$ is an $\perp$-type ;
\item[] - $X_{\perp}(t_1,...,t_n)$ is an $\perp$-type ;
\item[] - If $B$ is an $\perp$-type, then $A \f B$ is an $\perp$-type for every type $A$ ;
\item[] - If $A$ is an $\perp$-type, then $\q vA$ is an $\perp$-type for every variable $v$.
\end{itemize}

Therefore, $A$ is an $\perp$-type iff : $A = \q \sou{v_1} (E_1 \f F_1)$, $F_i=\q \sou{v_{i+1}}
(E_{i+1} \f F_{i+1})$ $1 \leq i \leq r-1$, and $F_r=\q \sou{v_{r+1}} X_{\perp}(t_1,...,t_n)$ or
$F_r=\q \sou{v_{r­1}}\perp$.\\

We add to the $AF2$ type system the new following rules :
\begin{itemize}  
\item[] (6$'$) If $\G\v t:A$, and $X_{\perp}$ has no free occurence in $\G$, then $\G\v t:\q
X_{\perp}A$.  
\item[] (7$'$) If $\G\v t:\q X_{\perp}A$, and $G$ is an $\perp$-type, then $\G\v t:A[G/X_{\perp}]$. 
\end{itemize} 
 
We call $AF2_{\perp}$ the new type system, and we write $\G\v_{\perp} t:A$ if $t$ is typable in
$AF2_{\perp}$ of type $A$ in the context $\G$.\\

\bf Remark \rm We can also see the system $AF2_{\perp}$ as a restriction of the system $AF2$.
Therefore, $AF2_{\perp}$ satisfies the same properties of $AF2$ (strongly normalization and
preservation of types). $\Box$

\subsection{The general Theorem}   

Let 
\begin{center}
$N^{\perp}[x] = \q X_{\perp} \{ X_{\perp}(0),\q y(X_{\perp}(y) \f X_{\perp}(sy)) \f
X_{\perp}(x) \}$
\end{center}  
where $X_{\perp}$ is a unary $\perp$-variable.

By the previous remark, we have : if  $\G\v_{\perp} t:N^{\perp}[s^n(0)]$, then $t \simeq\sb{\b}
\so{n}$.

\begin{lemma}
If $T$ is a closed normal $\l$-term such that $\v T: \q x \{N$*$[x] \f \neg\neg N[x]\}$, then
$\v_{\perp} T: \q x \{N^{\perp}[x] \f \neg\neg N[x] \}$.
\end{lemma}
\bf Proof \rm $T$ is a closed normal $\l$-term, then $T = \l\n T'$, and $ \n :N$*$[x] \v T':\neg\neg
N[x]$. Since $ \n :N^{\perp}[x] \v_{\perp} \n :N$*$[x] $, then $ \n :N^{\perp}[x] \v_{\perp}
T':\neg\neg N[x]$. Therefore  $\v_{\perp} T: \q x \{N^{\perp}[x] \f \neg\neg N[x] \}$.  $\Box$ \\
  
\bf Remarks \rm \\
1) We have $\v T_1,T_2: \q x \{N^{\perp}[x] \f \neg\neg N[x]\}$.\\
2) The $\l$-terms $T$ and $T'$ (given in the introduction) are of type  $\q x \{N^{\perp}[x] \f
\neg\neg N[x]\}$.
\begin{itemize}  
\item[] - We have $\n:N^{\perp}[x] \v_{\perp} \n:\perp,(\perp\f\perp)\f\perp$. Since $\n
:N^{\perp}[x] ,f:\neg N[x] \v_{\perp} (T_i) \n f: \perp$ and $\v_{\perp} \l xx:\perp \f \perp$, then
$\n :N^{\perp}[x] ,f:\neg N[x] \v_{\perp} ((\n)(T_i) \n f) \l xx:\perp$. Therefore $\v_{\perp} T:\q
x \{ N^{\perp}[x] \f \neg\neg N[x]\}$.
\item[] - We have $\n :N^{\perp}[x] \v_{\perp} \n: \perp ,( \perp \f \perp ) \f \perp$. Since  $\n
:N^{\perp}[x],f:\neg N[x] \v_{\perp} (T_i) \n f: \perp$ and $\n :N^{\perp}[x],f:\neg N[x] \v_{\perp}
\l d (T_i) \n f: \perp \f \perp$, then $\n :N^{\perp}[x],f:\neg N[x] \v_{\perp} ((\n)(T_i) \n f)\l
d(T_i)\n f:\perp$. Therefore $\v_{\perp} T':\q x\{ N^{\perp}[x] \f \neg\neg N[x]\}$. $\Box$
\end{itemize}

We give now a general type for storage operators for integers.

\begin{theo}
If $\v_{\perp} T: \q x \{ N^{\perp}[x] \f \neg\neg N[x] \}$, then $T$ is a storage operator for
the integers. 
\end{theo}

The type system $F_{\perp}$ is the subsystem of $AF2_{\perp}$ where we only have
propositional variables and constants (predicate variables or predicate symbols of arity 0). So, first
order variable, function symbols, and finite sets of equations are useless. The rules for typed
are 1), 2), 3), and 6), 7) restricted to propositional variables. For each predicate variable
(resp. predicate symbol) $X$, we associate a predicate variable (resp. a predicate symbol)
$X^{\di}$ of $F_{\perp}$ type system. For each formula $A$ of $AF2_{\perp}$, we associate the formula
$A^{\di}$ of $F_{\perp}$ obtained by forgetting in $A$ the first order part. If
$\G=x_1:A_1,...,x_n:A_n$ is a context of $AF2_{\perp}$, then we denote by $\G^{\di}$ the context
$x_1:A_1^{\di},...,x_n:A_n^{\di}$ of $F_{\perp}$.\\  We write $\G\v^{\di}_{\perp} t:A$ if $t$ is
typable in $F_{\perp}$ of type $A$ in the context $\G$.\\ 
We have obviously the following property : if $\G \v_{\perp} t:A$, then $\G^{\di}
\v^{\di}_{\perp} t:A^{\di}$.\\

Theorem 5.1 is a consequence of the following Theorem.

\begin{theo}
If $\v^{\di}_{\perp} T: N^{\perp} \f \neg\neg N$, then for every $n \geq 0$, there is an $m \geq
0$ and $\t_m \simeq\sb{\b} \so{m}$, such that for every  $\th_n \simeq\sb{\b} \so{n}$, there is a
substitution $\s$, such that $(T)\th_n f \p (f)\s(\t_m)$. 
\end{theo}

Indeed, if $\v_{\perp} T: \q x \{ N^{\perp}[x] \f \neg\neg N[x] \}$, then $\v^{\di}_{\perp} T:
N^{\perp} \f \neg\neg N$. Therefore for every $n \geq 0$, there is an $m \geq
0$ and $\t_m \simeq\sb{\b} \so{m}$, such
that for every  $\th_n \simeq\sb{\b} \so{n}$, there is a substitution $\s$, such that $(T)\th_n f
\p (f)\s(\t_m)$. We have $\v_{\perp} \so{n} : N[s^n(0)]$, then $ f:\neg N[s^n(0)] \v_{\perp} (T)
\so{n} f :\perp$, therefore $ f:\neg N[s^n(0)]\v_{\perp} (f)\so{m} :\perp$. By Theorem 2.2, we have
$\v_{\perp} \so{m} : N[s^n(0)]$ and thus $n=m$. Therefore $T$ is a storage operator for the
integers.  $\Box$ \\

In order to prove Theorem 5.2, we shall need some Lemmas.

\begin{lemma} If $\G,\n:N^{\perp}\v^{\di}_{\perp}(\n)\sou{d}:\perp$, then $\sou{d}=a,b,d_1,...,d_r$
and there is an $\perp$-type $F$, such that : $\G,\n:N^{\perp}\v^{\di}_{\perp} a:F$ ;
$\G,\n:N^{\perp}\v^{\di}_{\perp} b:F \f F$ ; $F \lhd E_1 \f F_1$, $F_i \lhd E_{i+1} \f F_{i+1}$
$1 \leq i \leq r-1$ ; $F_r \lhd \perp$ ; and $\G,\n:N^{\perp}\v^{\di}_{\perp} c_i : E_i$
$1 \leq i \leq r$.  
\end{lemma}
\bf Proof \rm We use Theorem 2.2.  $\Box$

\begin{lemma} If $F$ is an $\perp$-type and $\G,x:F \v^{\di}_{\perp} (x)\sou{d}:\perp$, then
$\sou{d}=d_1,...,d_r$ ; $F \lhd E_1 \f F_1$ ; $F_i \lhd E_{i+1} \f F_{i+1}$ $1 \leq i \leq r-1$ ; $F_r
\lhd \perp$ ; and $\G,x:F \v^{\di}_{\perp} c_i : E_i$ $1 \leq i \leq r$.  
\end{lemma} 
\bf Proof \rm We use Theorem 2.2.  $\Box$

\begin{lemma} Let $t$ be a normal $\l$-term, and $A_1,...,A_n$ a sequence of $\perp$-types.\\ 
If $x_1:A_1,...,x_n:A_n \v^{\di}_{\perp} t:N$, then there is an $m \geq 0$ such that $t = \so{m}$. 
\end{lemma}
\bf Proof \rm  We prove by induction on $u$ that if $u$ is a normal $\l$-term, $X$ a propositionnal
variable, and $x_1:A_1,...,x_n:A_n,x:X,f:X \f X \v^{\di}_{\perp} u:X$, then there is an $m \geq 0$
such that $u=(f)^m x$.  $\Box$ \\

We can now give the proof of Theorem 5.2.\\ 

\bf Proof of Theorem 5.2 \rm \\
Let $\n$ and $f$ two fixed variables, and $\v^{\di}_{\perp} T: N^{\perp} \f \neg\neg N$.\\
A good context $\G$ is a context of the form $\G = \n:N^{\perp}, f:\neg N, x_{n_1,a_1,b_1,\sou{c_1}} :
F_1 ,..., x_{n_p,a_p,b_p,\sou{c_p}} : F_p$ where $F_i$ is an $\perp$-type, and $\G \v^{\di}_{\perp}
a_i : F_i$, $\G \v^{\di}_{\perp} b_i : F_i \f F_i$, $0 \leq n_i \leq n-1$, and $1 \leq i \leq p$ .  \\

We will prove that for every $n \geq 0$, there is a finite sequence of head reduction $\{
U_i \p V_i \}_{1\leq i\leq r}$ such that :\\  
1) $U_1 = (T)\n f$ and $V_r = (f)\t$ where $\t \simeq\sb{\b}\so{m}$ for some $m \geq 0$ ;\\ 
2) $V_i = (\n)a b \sou{c}$ or $V_i = (x_{l,a,b,\sou{c}}) \sou{d}$ $0 \leq l \leq n-1$;\\ 
3) If $V_i = (\n)a b \sou{c}$, then $U_{i+1} = (a)\sou{c}$ if $n=0$ and
$U_{i+1} = ((b)x_{n-1,a,b,\sou{c}})\sou{c}$ if $n \neq 0$ \\  
4) If $V_i = (x_{l,a,b,\sou{c}})\sou{d}$ $0 \leq l \leq n-1$, then $U_{i+1}=(a)\sou{d}$ if $l=0$ and
$U_{i+1} = ((b)x_{l-1,a,b,\sou{d}})\sou{d}$ if $l \neq 0$.\\ 
5) There is a good context $\G$ such that $\G \v^{\di}_{\perp} V_i :\perp$ $1 \leq i \leq r$.\\ 

We have $\v^{\di}_{\perp} T: N^{\perp} \f \neg\neg N$, then  $\n:N^{\perp},f:\neg N
\v^{\di}_{\perp}(T)\n f:\perp$, and by Lemmas 5.2 and 5.3, $(T)\n f \p V_1$ where $V_1
= (f)\t$ or $V_1 = (\n)ab\sou{c}$.\\  
Assume that we have the head reduction $U_k \p V_k$ and $V_k \neq (f)\t$. 
\begin{itemize} 
\item[] 
- If $V_k = (\n)a b \sou{c}$, then, by induction hypothesis, there is a good context $\G$ such that
$\G \v^{\di}_{\perp}(\n)a b \sou{c} :\perp$. By Lemma 5.2, there is an $\perp$-type, such that $\G
\v^{\di}_{\perp}a:F$ ; $\G \v^{\di}_{\perp}b:F \f F$ ; $\sou{c}=c_1,...,c_s$ ; $F \lhd E_1 \f F_1$ ;
$F_i \lhd E_{i+1} \f F_{i+1}$ $1 \leq i \leq s-1$ ; $F_s \lhd \perp$ ; and $\G \v^{\di}_{\perp}c_i :
E_i$ $1 \leq i \leq s$.   
\begin{itemize}
\item[]
- If $n=0$, let  $U_{k+1} = (a)\sou{c}$. We have $\G \v^{\di}_{\perp} U_{k+1}:\perp$. 
\item[]
- If $n \neq 0$, let $U_{k+1} = ((b)x_{n-1,a,b,\sou{c}})\sou{c}$. The variable $x_{n-1,a,b,\sou{c}}$
is not used before. Indeed, if it is, we check easly that the $\l$-term $(T)\so{n} f$ is not solvable.
But that is impossible because $f:\neg N \v^{\di}_{\perp}(T)\so{n} f :\perp$. Let $\G' = \G
,x_{n-1,a,b,\sou{c}}:F$. $\G'$ is a good context and $\G' \v^{\di}_{\perp} U_{k+1} : \perp$. 
\end{itemize} 
\end{itemize}
\begin{itemize}
\item[] 
- If $V_k = (x_{l,a,b,\sou{c}}) \sou{d}$, then, by induction hypothesis, there is a good context $\G$
such that $\G \v^{\di}_{\perp}(x_{l,a,b,\sou{c}}) \sou{d} :\perp$.  $x_{l,a,b,\sou{c}} : F$ is in
the context $\G$, then by Lemma 5.3, $\sou{d}=d_1,...,d_s$ ; $F \lhd E_1 \f F_1$ ; $F_i \lhd E_{i+1}
\f F_{i+1}$ $1 \leq i \leq s-1$ ; $F_s \lhd \perp$ ; and $\G \v^{\di}_{\perp}d_i : E_i$
$1 \leq i \leq s$. 
\begin{itemize}
\item[] 
- If $l=0$, let $U_{k+1}=(a)\sou{d}$. We have $\G \v^{\di}_{\perp} U_{k+1}:\perp$.
\item[] 
- If $l \neq 0$,  let $U_{k+1} = ((b)x_{l-1,a,b,\sou{d}})\sou{d}$. The variable $x_{l-1,a,b,\sou{d}}$
is not used before. Indeed, if it is, we check easly that the $\l$-term $(T)\so{n} f$ is not solvable.
But that is impossible because $f:\neg N \v^{\di}_{\perp}(T)\so{n} f :\perp$. Let $\G' = \G
,x_{l-1,a,b,\sou{c}}:F$.  $\G'$ is a good context and $\G' \v^{\di}_{\perp} U_{k+1} :\perp$. 
\end{itemize}
\end{itemize}
Therefore there is a good context $\G'$ such that $\G' \v^{\di}_{\perp} U_{k+1} :\perp$, then, by
Lemmas 5.2 and 5.3, $U_{k+1} \p V_{k+1}$ where  $V_{k+1} = (f)\t$ or $V_{k+1} =
(\n)ab\sou{c}$ or  $V_{k+1} = (x_{l,a,b,\sou{c}})\sou{d}$ $0 \leq l \leq n-1$. \\ 
This constraction always terminates. Indeed, if not, we check easly that the $\l$-term $(T)\so{n} f$
is not solvable. But that is impossible because $f:\neg N \v^{\di}_{\perp}(T)\so{n} f :\perp$. \\ 
Therefore there is $r \geq 0$ and a good context $\G$ such that $\G \v^{\di}_{\perp} V_{r} = (f)\t
:\perp$, and by Theorem 2.2, $\G \v^{\di}_{\perp} \t :N$. Therefore by Lemma 5.4, there is an $m
\geq 0$ such that $\t \simeq\sb{\b}\so{m}$.\\
By the Theorem 4.2, we have the proof of the Theorem 5.2. $\Box$

\section{Generalization}

In this section, we give (without proof) a generalization of the Theorem 5.1.\\

Let $T$ be a closed $\l$-term, and $D,E$ two closed types of $AF2$ type system. We say that $T$ is a
storage operator for the pair of types $(D,E)$ iff for every $\l$-term $\v t:D$, there is $\l$-terms
$\t_t$ and $\t'_t$, such that $\t'_t \simeq\sb{\b} \t_t$, $\v\t_t:E$, and for every $\th_t
\simeq\sb{\b} t$, there is a substitution $\s$, such that $(T)\th_t f \p (f)\s(\t_t)$.\\

We define two sets of types of $AF2$ type system: $\O^+$ (set of $\q$-positive types), and $\O^-$ (set
of $\q$-negative types) in the following way :
\begin{itemize} 
\item[] - If $A$ is an atomic type, then $A \in \O^+$, and $A \in \O^-$ ;
\item[] - If $T \in \O^+$, and $T' \in \O^-$, then, $T' \f T \in \O^+$, and $T \f T' \in \O^-$ ; 
\item[] - If $T \in \O^+$, then $\q x T \in \O^+$ ;
\item[] - If $T \in \O^-$, then $\q x T \in \O^-$ ;
\item[] - If $T \in \O^+$, then $\q X T \in \O^+$ ;
\item[] - If $T \in \O^-$, and $X$ has no free occurence in $T$, then $\q X T \in \O^-$.
\end{itemize}
Therefore, $T$ is a $\q$-positive types iff the universal second order quantifier appears positively
in $T$. \\

For each predicate variable $X$, we associate an $\perp$- variable $X_{\perp}$.

For each formula $A$ of $AF2$ type system, we define the formula $A^{\perp}$ as follows : 
\begin{itemize} 
\item[] - If $A=\perp$, then  $A^{\perp} = A$ ;
\item[] - If $A=R(t_1,...,t_n)$, where $R$ is an $n$-ary predicate symbol, then  $A^{\perp} = A$ ;
\item[] - If $A=X(t_1,...,t_n)$, where $X$ is an $n$-ary predicate variable, then $A^{\perp} =
X_{\perp}(t_1,...,t_n)$; 
\item[] - If $A=B \f C$, then $A^{\perp} = B^{\perp} \f C^{\perp}$ ;
\item[] - If $A=\q x B$, then $A^{\perp} = \q x  B^{\perp}$ ;
\item[] - If $A=\q X B$, then $A^{\perp} = \q X_{\perp} B^{\perp}$.
\end{itemize}

$A^{\perp}$ is called the $\perp$-transformation of $A$.

\begin{theo} Let $D,E$ two $\q$-positive closed types of $AF2$ type system, such that $E$ does not
contain $\perp$. If $\v_{\perp} T: D^{\perp} \f \neg\neg E$, then $T$ is a storage operator for the
pair $(D,E)$.   
\end{theo}

\end{document}